\documentstyle[12pt]{article}
\parskip=\bigskipamount
\newtheorem{definition}{Definition}[section]
\newtheorem{lemma}[definition]{Lemma}
\newtheorem{theorem}[definition]{Theorem}

\newtheorem{remark}[definition]{Remark}
\newtheorem{example}[definition]{Example}

\font\fr=eufm10  scaled \magstep 1   
\font\ddpp=msbm10  scaled \magstep 1  

\def\QED{\hskip0.1em\hfill\null\ \null\nobreak\hfill
\kern3pt\lower1.8pt\vbox{\hrule\hbox   {\vrule\kern1pt\vbox{\kern1.7pt
\hbox{$\scriptstyle   QED$}\kern0.2pt}\kern1pt\vrule}\hrule}}
\def\R{\hbox{\ddpp R}}               

\begin{document}
\title{Characterization of Local Configuration Controllability for a 
Class of Mechanical Systems}         
\author{Jorge Cort\'es \and Sonia Mart{\'\i}nez}
\date{\today}          
\maketitle

\bigskip

\begin{abstract}
We investigate local configuration controllability for mechanical control systems within the affine connection formalism. Extending the work by Lewis for the single-input case, we are able to characterize local configuration controllability for systems with $n$ degrees of freedom and $n-1$ input forces.
\end{abstract}

\section{Introduction}
Mechanical control systems belong to a class of nonlinear systems whose controllability properties have not been fully characterized yet. Much work has been devoted to the study of their rich geometrical structure, both in the Hamiltonian framework (see \cite{NiSc} and references therein) and in the Lagrangian one, which is receiving increasing attention in the last years \cite{BuLe,KeMu,Le,Le2,LeMu2,LeMu3,OsBu}. This research is providing new insights and a bigger understanding of the accessibility and controllability aspects associated to them. In this direction, the affine connection formalism has revealed to be very useful modelling different types of mechanical systems, such as natural ones (Lagrangian equal to kinetic energy minus potential energy) \cite{LeMu2,LeMu3}, with symmetries \cite{BuLe}, with nonholonomic constraints \cite{Le2},... and, on the other hand, it has led to the development of some new techniques and control algorithms for approximate trajectory generation in controller design \cite{BuLeLe,BuLeo}. Certainly, we shall see further progress in these directions in future years.

From a theoretical point of view, underactuated mechanical control systems offer a control challenge as they have non-zero drift, their linearisation at zero velocity is not controllable and they are not feedback linearisable  as well. That is, they are not amenable to standard techniques in control theory \cite{NiSc}. The work by Lewis and Murray \cite{LeMu2,LeMu3} on simple mechanical control systems has rendered strong conditions for configuration accessibility and sufficient conditions for configuration controllability. The conditions for the latter one are based on the sufficient conditions that Sussmann obtained for general affine control systems \cite{Su2}. It is worthy to note the fact that these conditions are not invariant under input transformations, as the simple example of the planar rigid body shows \cite{LeMu2}.

As the controllability is the more interesting property in practice, more research is needed in order to sharpen the configuration controllability conditions. Lewis \cite{Le} started this process by investigating and fully solving the single-input case, building on previous results by Sussmann for general scalar-input systems \cite{Su}. The recent work by Bullo \cite{Bu,Bu2} on series expansion for the evolution of a mechanical control system has given the necessary tools to tackle this problem in the much more involved multi-input case, as well as a powerful machinery for the design of motion control algorithms for a large class of autonomous vehicles, robotic manipulators and locomotion devices. In this paper, we characterize local configuration controllability when the number of inputs and the degrees of freedom differs by one. Both results, Lewis' and ours, can be seen as particular cases of the following conjecture, which remains open: {\sl The system is local configuration controllable at a point if and only if there exists a basis of inputs satisfying the sufficient conditions for local configuration controllability at that point}. The conjecture relies on the fact we have mentioned before: the lack of invariance of the sufficient conditions under input transformations. It is remarkable to note that local controllability has not been characterized yet for general control systems, even for the single input case (in this respect see \cite{He,Su,Su2}).

The paper is organised as follows. In Section 2, we describe the affine connection framework for mechanical control systems and recall the controllability notions we shall consider on them. In Section 3 we review the existing results concerning configuration controllability \cite{LeMu2,LeMu3} and the series expansion for the evolution of a mechanical control system developed by Bullo in \cite{Bu,Bu2}. In Section 4 we briefly recall the single-input case solved by Lewis and properly state his conjecture. Section 5 contains the main contributions of this paper. Finally, in Section 6 we expose a simple example on a $4$-dimensional manifold with $2$ inputs whose we can not decide if it is configuration controllable or not. It certainly does not admit a new input basis verifying the sufficient conditions, so that if controllable, it would be a counter example of the above stated conjecture.

\section{Simple Mechanical Control Systems}

Let $Q$ be a $n$-dimensional manifold. We will denote by $TQ$ the tangent bundle of $Q$, by $\hbox{\fr X} (Q)$ the set of vector fields on $Q$ and by $C^{\infty}(Q)$ the set of smooth functions on $Q$.

A {\bf simple mechanical control system} is defined by a triple $(Q,g,{\cal F})$, where $Q$ is the manifold of configurations of the system, $g$ is a Riemannian metric on $Q$ and ${\cal F}=\{F^1,\dots,F^m\}$ is a set of $m$ linearly independent $1$-forms on $Q$, which physically correspond to forces or torques.

Associated with the Riemannian metric $g$ there is a natural affine connection, called the Levi-Civita connection. An {\bf affine connection} \cite{AbMa,KoNo} is defined as an assignment
\[
\begin{array}{rccc}
\nabla: & \hbox{\fr X}(Q) \times \hbox{\fr X} (Q) & \longrightarrow & \hbox{\fr X} (Q) \\
& (X , Y) & \longmapsto & \nabla_X Y
\end{array}
\]
which is $\R$-bilinear and satisfies $\nabla_{fX} Y =f \nabla_X Y$ and $\nabla_X(fY)=f \nabla_XY + X(f)Y$, for any $X$, $Y \in \hbox{\fr X}(Q)$, $f\in C^{\infty}(Q)$. A curve $c:[a,b] \longrightarrow Q$ is a {\bf geodesic} for $\nabla$ if $\nabla_{\dot{c}(t)}\dot{c}(t)=0$. Locally, the condition for a curve $t \rightarrow (x^1(t),\dots,x^n(t))$ to be a geodesic can be expressed as
\begin{equation}\label{geodesic}
\ddot{x}^a + \Gamma^a_{bc} \dot{x}^b \dot{x}^c =0 \, , \; \; 1 \le a \le n \, ,
\end{equation}
where the $\Gamma^a_{bc}(x)$ are the Christoffel symbols of the affine connection, that is, they are given by $\displaystyle{\nabla_{\frac{\partial}{\partial x^b}} \frac{\partial}{\partial x^c} = \Gamma^a_{bc} \frac{\partial}{\partial x^a}}$. The geodesic equation (\ref{geodesic}) is a first-order differential equation on $TQ$. The vector field corresponding to this first-order equation is given in coordinates by
\[
S=v^a \frac{\partial}{\partial x^a} - \Gamma^a_{bc} v^b v^c\frac{\partial}{\partial v^a} \, .
\]
and is called the {\bf geodesic spray} associated with the affine connection $\nabla$. Hence, the integral curves of the geodesic spray $S$, $(x^a,\dot{x}^a)$ are the solutions of the geodesic equation.

The Levi-Civita connection $\nabla^g$ is determined by the formula
\begin{eqnarray*}
g(\nabla^g_XY,Z) &=& \frac{1}{2} \left( X(g(Y,Z)) +Y(g(Z,X)) -Z(g(X,Y)) \right. \\
&& \left. + g(Y,[Z,X]) - g(X,[Y,Z]) + g(Z,[X,Y]) \right) \, , \; \; X,Y,Z \in \hbox{\fr X}(Q) \, .
\end{eqnarray*}

Instead of the input forces $F^1,\dots,F^m$, we shall make use of the vector fields $Y_1,\dots, Y_m$, defined as $Y_i=\sharp_g (F^i)$, where $\sharp_g = \flat_g^{-1}$ and $\flat_g:TQ \longrightarrow T^*Q$ is the musical isomorphism given by $\flat_g(X)(Y)=g(X,Y)$. Roughly speaking, this corresponds to consider ``accelerations" rather than forces.

The control equations for the simple mechanical control system may then be written as
\begin{equation}\label{affine}
\nabla^g_{\dot{c}(t)} \dot{c}(t) = \sum u^i(t) Y_i(c(t)) \, .
\end{equation}
The inputs we consider come from the set
\[
{\cal U} = \{ u:[0,T] \rightarrow \R^m | \, T>0, u \; \hbox{is measurable and} \, \|u\| \le 1 \} \, .
\]
We can use a general affine connection in (\ref{affine}) instead of the Levi-Civita connection without changing the structure of the equation. This is particularly interesting, since nonholonomic mechanical control systems give also rise to equations of the form (\ref{affine}) as explained in \cite{Le2}.

We can turn (\ref{affine}) into a general affine control system with drift
\begin{equation}\label{drift}
\dot{x}(t) = f(x(t)) + \sum u^i(t) g_i(x(t)).
\end{equation}
To do this we need another bit of notation. The vertical lift of a vector field $X$ on $Q$ is the vector field $X^v$ on $TQ$ defined as
\[
X^v(v_q) = \frac{d}{dt}\Big|_{t=0}(v_q+ t X(q)) \, .
\]
In coordinates, if $\displaystyle{X=X^a \frac{\partial}{\partial q^a}}$, one can check that $\displaystyle{X^v=X^a \frac{\partial}{\partial v^a}}$. Then, the second-order equation (\ref{affine}) on $Q$ can be written as the first-order system on $TQ$
\begin{equation}\label{conver}
\dot{v} = S(v) + \sum u^i(t) Y_i^v(v) \, ,
\end{equation}
where $S$ is the geodesic spray associated with the affine connection $\nabla^g$.

\subsection{Controllability Notions}

The control equations for the mechanical system (\ref{conver}) are nonlinear. No standard technique in control theory \cite{NiSc}, as for example the linearisation around an equilibrium point or linearisation by feedback, yields satisfactory results in the analysis of its controllability properties.

The point in the approach of Lewis and Murray to simple mechanical control systems is precisely to know what is happening to configurations, rather than to states, since in many of these systems, configurations may be controlled, but not configurations and velocities at the same time. The basic question they pose is ``what is the set of configurations which are attainable from a given configuration starting from rest?" Moreover, since we are dealing with objects defined on the configuration manifold $Q$, we expect to find answers on $Q$, although the control system (\ref{conver}) lives in $TQ$.

\newpage

\begin{definition}
A {\bf solution} of (\ref{affine}) is a pair $(c,u)$, where $c:[0,T] \longrightarrow Q$ is a piecewise smooth curve and $u\in {\cal U}$ such that $(\dot{c},u)$ satisfies the first order control system (\ref{conver}).
\end{definition}

Consider $q_0 \in Q$, $(q_0,0_{q_0}) \in T_{q_0}Q$ and let $U \subset Q$, $\bar{U} \subset TQ$ be neighbourhoods of $q_0$ and $(q_0,0_{q_0})$, respectively. Define
\begin{eqnarray*}
{\cal R}_Q^U (q_0,T) = \left\{ q \in Q \hspace{-3pt} \right. &|& \left. \hspace{-3pt} \hbox{there exists a solution $(c,u)$ of (\ref{affine}) such that} \right. \\
&&  \left. \dot{c}(0)=0_{q_0} , c(t) \in U \, \hbox{for} \; t \in [0,T] \,\hbox{and} \; \dot{c}(T) \in T_qQ \right\} \\
{\cal R}_{TQ}^{\bar{U}} (q_0,T) = \left\{ (q,v) \in TQ \hspace{-3pt} \right. &|& \left. \hspace{-3pt} \hbox{there exists a solution $(c,u)$ of (\ref{affine}) such that} \right. \\
&&  \left. \dot{c}(0)=0_{q_0} , (c(t),\dot{c}(t)) \in \bar{U} \, \hbox{for} \; t \in [0,T] \,\hbox{and} \; \dot{c}(T)=v \in T_qQ \right\}
\end{eqnarray*}

and denote by ${\cal R}_Q^U (q_0, \le T) = \cup_{0 \le t \le T} {\cal R}_Q^U (q_0,t)$, ${\cal R}_{TQ}^{\bar{U}} (q_0, \le T) = \cup_{0 \le t \le T} {\cal R}_{TQ}^{\bar{U}} (q_0,t)$.

Now, we recall the notions of accessibility considered in \cite{LeMu2}.

\begin{definition}
The system (\ref{affine}) is {\bf locally configuration accessible (LCA) at $q_0 \in Q$} if there exists $T>0$ such that ${\cal R}_Q^U (q_0, \le t)$ contains a non-empty open set of $Q$, for all neighbourhoods $U$ of $q_0$ and all $0\le t \le T$. If this holds for any $q_0 \in Q$ then the system is called locally configuration accessible.
\end{definition}

\begin{definition}
The system (\ref{affine}) is {\bf locally accessible (LA) at $q_0 \in Q$ and zero velocity} if there exists $T>0$ such that ${\cal R}_{TQ}^{\bar{U}} (q_0, \le t)$ contains a non-empty open set of $TQ$, for all neighbourhoods $\bar{U}$ of $(q_0,0_{q_0})$ and all $0\le t \le T$. If this holds for any $q_0 \in Q$ then the system is called locally accessible at zero velocity.
\end{definition}

We shall focus our attention on the following concepts of controllability \cite{LeMu2}.

\begin{definition}
The system (\ref{affine}) is {\bf small-time locally configuration controllable (STLCC) at $q_0 \in Q$} if there exists $T>0$ such that ${\cal R}_Q^U (q_0, \le t)$ contains a non-empty open set of $Q$ to which $q_0$ belongs, for all neighbourhoods $U$ of $q_0$ and all $0\le t \le T$. If this holds for any $q_0 \in Q$ then the system is called small-time locally configuration controllable.
\end{definition}

\begin{definition}
The system (\ref{affine}) is {\bf small-time locally controllable (STLC) at $q_0 \in Q$ and zero velocity} if there exists $T>0$ such that ${\cal R}_{TQ}^{\bar{U}} (q_0, \le t)$ contains a non-empty open set of $TQ$ to which $(q_0,0_{q_0})$ belongs, for all neighbourhoods $\bar{U}$ of $(q_0,0_{q_0})$ and all $0\le t \le T$. If this holds for any $q_0 \in Q$ then the system is called small-time locally controllable at zero velocity.
\end{definition}

\section{Existing Results}

Here we review some accessibility and controllability results obtained in \cite{LeMu2,LeMu3} and expose the work by Bullo \cite{Bu,Bu2} in describing the evolution of mechanical control systems via a series expansion. This series will be key in the proof of the main results of this paper. 

\subsection{On Controllability}

The {\bf symmetric product} on $\hbox{\fr X}(Q)$ is defined by
\[
\left \langle X:Y\right \rangle = \nabla_XY + \nabla_YX \, .
\]
The geometric meaning of the symmetric product is the following \cite{Le3}: a {\bf geodesically invariant} distribution ${\cal D}$ is a distribution such that for every geodesic
$c(t)$ of $\nabla$ starting from a point in ${\cal D}$, $\dot{c}(0) \in {\cal D}_{c(0)}$, we have that $\dot{c}(t) \in {\cal D}_{c(t)}$. Then, one can prove that ${\cal D}$ is geodesically invariant if and only if $\langle X:Y \rangle \in {\cal D}$, $\forall X$, $Y \in {\cal D}$. 

Given the input vector fields ${\cal Y}=\{ Y_1,\dots , Y_m \}$, let us denote by $\overline{Sym}(\cal Y)$ the distribution obtained by closing the set ${\cal Y}$ under the symmetric product and by $\overline{Lie}(\cal Y)$ the involutive closure of ${\cal Y}$. With these ingredients, one can prove

\begin{theorem} (\cite{LeMu2})
The control system (\ref{affine}) is locally configuration accessible at $q$ (respectively locally accessible at $q$ and zero velocity) if $\overline{Lie}(\overline{Sym}({\cal Y}))_q=T_qQ$ (respec. $\overline{Sym}({\cal Y})_q=T_qQ$).
\end{theorem}

If $P$ is a symmetric product of vector fields in ${\cal Y}$, we let $\gamma_i(P)$ denote the number of occurrences of $Y_i$ in $P$. The {\bf degree} of $P$ will be $\gamma_1(P) + \dots + \gamma_m(P)$. We shall say that $P$ is {\bf bad} if $\gamma_i(P)$ is even for each $1 \le i \le m$. We say that $P$ is {\bf good} if it is not bad. The following theorem gives sufficient conditions for STLCC.

\begin{theorem}\label{suffi}
Suppose that the system is LCA at $q$ (respectively, LA at $q$ and zero velocity) and that ${\cal Y}$ is such that every bad symmetric product $P$ at $q$ in ${\cal Y}$ can be written as a linear combination of good symmetric products at $q$ of lower degree than $P$. Then (\ref{affine}) is STLCC at $q$ (respec. STLC at $q$ and zero velocity).
\end{theorem}

This theorem was proved in \cite{LeMu2}, adapting previous work by Sussmann \cite{Su2} on general control systems of the form (\ref{drift}).

\subsection{Series Expansion}

We would like to have a description of the evolution of the mechanical control system (\ref{affine}) when starting from rest. This is accomplished in \cite{Bu,Bu2} via a series expansion on the configuration manifold $Q$.
 
In the sequel, we let
\[
Z(q,t)=\sum_{i=1}^{m} u_{i}(t)Y_i(q) \, .
\]
We have the following
\begin{theorem}
Let $c(t)$ be the solution of equation (\ref{affine}) with input given by $Z(q,t)$ and with initial conditions $c(0)=q_0$, $\dot{c}(0)=0$. Let the Christoffel symbols $\Gamma_{jk}^i (q)$ and the vector field $Z(q,t)$ be uniformly integrable and bounded analytic in a neighbourhood of $q_0$. Define recursively the time varying vector fields 
\begin{eqnarray*}
V_1(q,t) &=& \int_0^t Z(q,s)ds \, , \\
V_k (q,t) &=& -\frac{1}{2} \sum_{j=1}^{k-1} \int_0^t \left \langle V_j (q,s): V_{k-j}(q,s) \right \rangle ds \, , \; k \geq 2\, ,
\end{eqnarray*}
where $q$ is mantained fixed at each integral. Then there exists a sufficiently small $T_c$ such that the series $\sum_{k=1}^{\infty}V_{k}(q,t)$ converges absolutely and uniformly in $t$ and $q$, for all $t \in [0,T_c]$ and for all $q$ in a neighbourhood of $q_0$. Over the same interval, the solution $c(t)$ satisfies 
\begin{equation}\label{bullo}
\dot{c}(t) = \sum_{k=1}^{\infty} V_k(c(t),t) \, .
\end{equation}
\end{theorem}

This theorem generalizes various previous results obtained in \cite{BuLeLe,BuLeo} under the assumption of small amplitude forcing. The first few terms of the series (\ref{bullo}) can be computed to obtain 
\begin{eqnarray}\label{series}
\dot{c}(t) \hspace{-3pt} &=& \hspace{-3pt} \overline{Z}(c(t),t)- \frac{1}{2} \overline{ \left \langle \overline{Z}: \overline{Z}\right \rangle} (c(t),t) + \frac{1}{2} \overline{\left \langle \overline{\left \langle\overline{Z}:\overline{Z}\right \rangle} : \overline{Z}\right \rangle} (c(t),t) \nonumber \\
\hspace{-4pt}&&\hspace{-4pt} - \frac{1}{2} \overline{\left \langle\overline{\left \langle\overline{\left \langle\overline{Z}:\overline{Z}\right \rangle}:\overline{Z}\right \rangle}:\overline{Z}\right \rangle} (c(t),t) -  \frac{1}{8} \overline{\left \langle\overline{\left \langle\overline{Z}:\overline{Z}\right \rangle}:\overline{\left \langle\overline{Z}:\overline{Z}\right \rangle}\right \rangle} (c(t),t) + O(\|Z\|^5 t^9) , 
\end{eqnarray}
where $\overline{Z}(q,t) \equiv \int_0^t Z(q,s)ds$ and so on.

\section{The One-Input Case}

Theorem \ref{suffi} gives us sufficient conditions for small-time local configuration controllability. A natural concern both from the theoretical and the practical point of view is to treat to sharpen this controllability test. Lewis \cite{Le} investigated the single-input case and managed to prove the next

\begin{theorem}\label{Lewis}
Let $(Q,g)$ be an analytic manifold with an affine connection $\nabla$. Let $Y$ be an analytic vector field on $Q$ and $q_0 \in Q$. Then the system
\[
\nabla_{\dot{c}(t)} \dot{c}(t)=u(t)Y(c(t))
\]
is locally configuration controllable at $q_0 \in Q$ if and only if $\hbox{dim}\,Q=1$.
\end{theorem}

The fact of being able to completely characterize STLCC in the single-input case (something which has not been accomplished yet for general control systems of the form (\ref{drift})) suggests that understanding local configuration controllability for mechanical systems may be possible. More precisely, examining the single-input case, one can deduce that if (\ref{affine}) is STLCC at $q_0$ then $\hbox{dim}\,Q=1$, which implies $\left \langle Y:Y \right \rangle(q_0) \in \hbox{span}\{Y(q_0)\}$, i.e. sufficient conditions for STLCC are also necessary. Can this be extrapolated to the multi-input case? The following conjecture was posed by Lewis:
\begin{center}
{\sl Consider a mechanical control system of the form (\ref{affine}) which is locally configuration accessible at $q_0 \in Q$. Then it is STLCC at $q_0$ if and only if there exists a basis of input vector fields which satisfies the sufficient conditions for STLCC at $q_0$.
}
\end{center}

Theorem \ref{Lewis} implies that the conjecture is true for $m=1$. In the following section we prove that this conjecture is also valid for $m=n-1$.

\section{The case $m=n-1$}

The following lemma, taken from \cite{Su}, will be helpful in the proof of the theorem of this section.

\begin{lemma}\label{magic}
Let $Q$ be an analytic manifold. Given $q_0 \in Q$ and $X_1,\dots,X_p \in \hbox{\fr X}(Q)$, 
$p\le n$, linearly independent vector fields, there exists a function $\phi:Q \longrightarrow 
\R$ satisfying the properties
\begin{enumerate}
\item $\phi$ is analytic
\item $\phi(q_0)=0$
\item $X_1(\phi)=\dots=X_{p-1}(\phi)=0$ on a neighbourhood $V$ of $q_0$
\item $X_p(\phi)(q_0)=-1$
\item Within any neighbourhood of $q_0$ there exists points $q$ where $\phi(q) < 0$ and $\phi(q)>0$.
\end{enumerate}
\end{lemma}

{\bf Proof:}
Let $Z_1,\dots,Z_n$ be vector fields in a neighbourhood of $q_0$ such that $\{Z_1(q_0),\dots,Z_n(q_0)\}$ form a basis for $T_{q_0}Q$ and $Z_i=X_i$, $1 \le i \le p-1$, $Z_p=-X_p$. Let $t_i \longmapsto \Psi_i(t)$ be the flow of $Z_i$, $1 \le i \le n$. In a sufficiently small neighbourhood $V$ of $q_0$, any point $q$ mat be expressed as $q=\Psi_1(t_1)\circ \dots \circ \Psi_n(t_n)(q_0)$ for some unique $n$-tuple $(t_1,\dots,t_n)\in \R^n$. Define $\phi(q)=t_p$. It is a simple exercise to verify that $\phi$ satisfies the required properties.
\QED

Next, we prove the main result of this paper. To get an idea of how the proof works, the interested reader is referred to \cite{CoMa}, where the simpler case $m=n-1=2$ was treated. The proof of the general result is considerable more involved that the one presented there (due to the higher number of inputs) and consists of a careful refinement of it.

\begin{theorem}\label{ole}
Consider a mechanical control system on a $n$-dimensional configuration
manifold $Q$ of the form (\ref{affine}) with $n-1$ inputs, which is locally
configuration accessible at $q_0 \in Q$. Then the system is locally
configuration controllable at $q_0$ if and only if there exists a basis
input vector fields satisfying the sufficient conditions for STLCC at $q_0$.
\end{theorem}
{\bf Proof:} There are two possibilities
\begin{itemize}
\item{ $ \forall Y_1, Y_2 \in D$, $\left \langle Y_1:
Y_2 \right \rangle(q_0)\in D_{q_0}$. }
\item There exist $Y_1$, $Y_2 \in D$
such that $\left \langle Y_1: Y_2\right \rangle (q_0) \not \in D_{q_0}$.
\end{itemize}
In the first case, there is nothing to prove. In the second one, we distinguish again two possibilities
\begin{enumerate}
\item{ There exists $Y_1 \in D$ with $\left \langle Y_1: Y_1\right \rangle(q_0) \not \in D_{q_0}$.}
\item There exist $Y_1$, $Y_2 \in D$, linearly independent at $q_0$ and such
that $\left \langle Y_1: Y_2\right \rangle (q_0) \not \in D_{q_0}$.
\end{enumerate}
Case (i) can be easily reduced to case (ii): if $\left \langle Y_1: Y_2\right \rangle (q_0) 
\in D_{q_0}$, then define a new $Y_2'$ by $Y_1+Y_2$ and we are done. So
let us treat the second case. We can complete the set $\{ Y_1(q_0), Y_2(q_0) \}$ to a basis
of $D_{q_0}$, $\{Y_1(q_0), Y_2(q_0), \dots, Y_m(q_0) \}$. Then, we have that $\hbox{span} \{Y_1(q_0), Y_2(q_0), \dots, Y_m(q_0), \left \langle Y_1:Y_2 \right \rangle (q_0) \} =
T_{q_0}Q$ and can write
\begin{eqnarray*}
\left \langle Y_1:Y_1 \right \rangle
(q_0) &=& lc (Y_1(q_0),\dots, Y_m(q_0)) + a_{11} \left \langle Y_1:Y_2
\right \rangle (q_0) \\
& \vdots & \\
\left \langle Y_m:Y_m \right \rangle
(q_0) &=& lc (Y_1(q_0),\dots, Y_m(q_0)) + a_{mm} \left \langle Y_1:Y_2
\right \rangle (q_0) \\
\left \langle Y_1:Y_3 \right \rangle (q_0) &=& lc
(Y_1(q_0),\dots, Y_m(q_0)) + a_{13}\left \langle Y_1:Y_2 \right \rangle
(q_0) \\
& \vdots & \\
\left \langle Y_{m-1}:Y_m \right \rangle (q_0) &=& lc
(Y_1(q_0),\dots, Y_m(q_0)) + a_{m-1m}\left \langle Y_1:Y_2 \right \rangle
(q_0) \, ,
\end{eqnarray*}
where $lc(Y_1(q_0),\dots,Y_m(q_0))$ means a linear
combination of $Y_1(q_0),\dots,Y_m(q_0)$. If $a_{11}=\dots=a_{mm}=0$, we have finished. Suppose then that there exists $s=s_1$ such that $a_{s_1s_1} \not = 0$. To prove the theorem, we have to find a change of basis $B=(b_{jk})$, $\hbox{det} \, (b_{jk}) \not = 0$, providing new vector fields in $D$
\[
Y_j' = \sum_{k=1}^{m} b_{jk} Y_k \, , \; \; 1 \le j \le m \, ,
\]
satisfying the sufficient conditions for STLCC at $q_0$. Since
\begin{eqnarray}\label{eluno}
\left \langle Y_j':Y_j' \right \rangle (q_0)&=&
\sum_{k,l=1}^{m} b_{jk} b_{jl} \left \langle Y_k:Y_l \right \rangle
(q_0)=\sum_{k=1}^{m} b_{jk}^2 \left \langle Y_k:Y_k \right \rangle (q_0) + 2
\hspace{-8pt} \sum_{1 \le k < l \le m} b_{jk} b_{jl} \left \langle Y_k:Y_l
\right \rangle (q_0) \nonumber \\
&=& lc(Y_1'(q_0),\dots,Y_m'(q_0)) + \left(
\sum_{k=1}^{m} b_{jk}^2 a_{kk} + 2 \hspace{-8pt} \sum_{1 \le k < l \le m}
b_{jk} b_{jl} a_{kl} \right) \left \langle Y_1:Y_2 \right \rangle (q_0)\, ,
\end{eqnarray}
the matrix $B$ must fulfill
\begin{equation}\label{itera1}
\sum_{k=1}^{m} b_{jk}^2 a_{kk} + 2 \hspace{-8pt} \sum_{1 \le k < l \le m} b_{jk} b_{jl} a_{kl} = 0 \, , \; \; 1 \le j \le m \, ,
\end{equation}
which is equivalent to
\begin{eqnarray*}
b_{js_1}&=&\displaystyle{\frac{-\sum_{k\not =s_1} b_{jk} a_{ks_1}}{a_{s_1s_1}}} \\
&& \pm \displaystyle{\frac{\sqrt{ (\sum_{k\not =s_1} b_{jk} a_{ks_1})^2 - a_{s_1s_1} (\sum_{k\not = s_1}b_{jk}^2 a_{kk}+2 \sum_{k < l, k,l \not= s_1} b_{jk} b_{jl} a_{kl})}}{a_{s_1s_1}}} \, , \; \; 1 \le j \le m \, ,
\end{eqnarray*}
due to $a_{s_1s_1} \not = 0$. After some computations, the radicand of this expression becomes 
\begin{eqnarray*}
\sum_{k\not =s_1} \hspace{-2pt} b_{jk}^2 (a^2_{ks_1} -a_{s_1s_1} a_{kk})
+ 2 \hspace{-9pt} \sum_{k<l, k,l \not = s_1} \hspace{-8pt} b_{jk} b_{jl} (a_{ks_1}a_{ls_1} 
-a_{s_1s_1}a_{kl}) \, .
\end{eqnarray*}
Denoting by
\begin{eqnarray*}
a^{(2)}_{kl} = a_{ks_1} a_{ls_1} -a_{s_1s_1}a_{kl} \, , \; \; k,l \in \{ 1, \dots, m \} / \{ s_1\} \, ,
\end{eqnarray*}
we have that the radicand would vanish if
\begin{equation}\label{itera2}
\sum_{k \not= s_1} b_{jk}^2 a_{kk}^{(2)} + 2 \hspace{-8pt} \sum_{k < l, k,l \not =s_1} 
b_{jk} b_{jl} a_{kl}^{(2)} = 0 \, .
\end{equation}
Note the similarity between (\ref{itera1}) and (\ref{itera2}). Now, several situations can occur
that we study in the following.

Suppose that there exists a $s_2$ such that $a_{s_2s_2}^{(2)} \not = 0$. Then we can repeat the same steps. Define recursively
\begin{equation}\label{defin}
\begin{array}{rcl}
a_{kl}^{(1)} & = & a_{kl} \, , \\
a^{(i)}_{kl} & = & a^{(i-1)}_{ks_{i-1}}a^{(i-1)}_{ls_{i-1}} - a^{(i-1)}_{s_{i-1}s_{i-1}} 
a^{(i-1)}_{kl} \, , \; i \ge 2 \, , \; \; k, l \in \{ 1, \dots, m \}/ \{ s_1, \dots, s_{i-1} \} \, .
\end{array}
\end{equation}
Reasoning as before, we obtain that (\ref{itera2}) would imply that
\[ 
b_{js_2}= lc (b_{j1},\dots, \hat{b}_{js_1}, \dots, \hat{b}_{js_2}, \dots, b_{jm}) \pm
\frac{1}{a_{s_2s_2}^{(2)}} \sqrt{\sum_{k \not = s_1,s_2} b_{jk}^2 a^{(3)}_{kk} + 2
\hspace{-8pt} \sum_{k < l,k,l \not =s_1, s_2} b_{jk} b_{jl} a_{kl}^{(3)} } \, , \;
1 \le j \le m \, ,
\]
where the symbol $\hat{b}$ means that that the term $b$ has been removed. This procedure can be 
iterated, assuming always that there exists $s_i$ such that $a^{(i)}_{s_is_i} \not = 0$. In this way, we finally obtain the following equations for the $b_{js_{m-1}}$,
\[ 
b_{js_{m-1}}=b_{js_{m}} \frac{-
a_{s_{m-1}s_{m}}^{(m-1)} \pm \sqrt{(a_{s_{m-1}s_{m}}^{(m-1)})^2 - a_{s_{m-1}s_{m-1}}^{(m-1)}
a_{s_ms_{m}}^{(m-1)}}}{a_{s_{m-1}s_{m-1}}^{(m-1)}} \, , \; \; 1 \le j \le m \, . 
\]
Let $(b_{js_m})_{1 \le j \le m}$ be a non-zero vector in $\R^m$. Now, if $(a_{s_{m-1}s_m}^{(m-1)})^2 - a_{s_{m-1}s_{m-1}}^{(m-1)} a_{s_ms_m}^{(m-1)} >0$, the quadratic polynomial in $b_{js_{m-1}}$
\begin{equation}\label{leches}
a_{s_{m-1}s_{m-1}}^{(m-1)}b_{js_{m-1}}^2 + 2a_{s_{m-1}s_m}^{(m-1)} b_{js_{m-1}}b_{js_m} + a_{s_ms_m}^{(m-1)}b_{js_m}^2 \, ,
\end{equation}
has two real roots and we can choose $(b_{js_{m-1}})_{1 \le j \le m} \in \R^m$ linearly independent with $(b_{js_m})_{1 \le j \le m}$ and such that (\ref{leches}) be positive for all $1 \le j \le m$. As this polynomial is the radicand of the preceding one,
\begin{equation}\label{next}
\sum_{k \not = s_1,\dots,s_{m-3}} \hspace{-8pt} b_{jk}^2 a_{kk}^{(m-2)} + 2 \hspace*{-23pt} \sum_{k<l,k,l \not = s_1,\dots,s_{m-3}} \hspace*{-17pt} b_{jk} b_{jl} a_{kl}^{(m-2)} \, ,
\end{equation}
our choice of $(b_{js_{m-1}})_{1 \le j \le m}$ ensures that we can again take $(b_{js_{m-2}})_{1 \le j \le m} \in \R^m$, linearly independent with $(b_{js_{m-1}})_{1 \le j \le m}$ and $(b_{js_m})_{1 \le j \le m}$ such that (\ref{next}) is positive for all $1 \le j \le m$. This is inherited step by step through the iteration process and we are able to choose a non-singular matrix $(b_{jk})$ satisfying (\ref{itera1}).

If $(a_{s_{m-1}s_m}^{(m-1)})^2 - a_{s_{m-1}s_{m-1}}^{(m-1)} a_{s_ms_m}^{(m-1)} <0$, then (\ref{leches}) does not change its sign for all $b_{js_{m-1}}$, $b_{js_m}$. Indeed, we have that $\hbox{sign} (\ref{leches}) = \hbox{sign} (a_{s_{m-1}s_{m-1}}^{(m-1)} ) = \hbox{sign} (a_{s_ms_m}^{(m-1)})$, $1 \le j \le m$. If this sign is positive, the former argument ensures us the choice of the desired matrix. If negative, it implies that (\ref{next}) does not change its sign for all $b_{js_{m-2}}$, $b_{js_{m-1}}$, $b_{js_m}$. In particular, notice that this ensures us that $\hbox{sign} (\ref{next}) = \hbox{sign} (a_{s_{m-2}s_{m-2}}^{(m-2)} ) = \hbox{sign} (a_{s_{m-1}s_{m-1}}^{(m-2)}) = \hbox{sign} (a_{s_ms_m}^{(m-2)})$, for all $1 \le j \le m$. Then, the unique problem we must face is when, through the iteration process, all the radicands are negative. In this case, we can apply Lemma \ref{magic} to the vector fields $\{Y_1,\dots,Y_m,\langle Y_1:Y_2 \rangle \}$ to find a 
function $\phi$ satisfying the properties {\it (i)-(v)}. By (\ref{series}), we have 
that
\begin{eqnarray*}
\dot{c}(t) &=& \sum_{i=1}^m \bar{u}_i Y_i - \frac{1}{2} \overline{\langle 
\sum_{j=1}^m \bar{u}_j Y_j : \sum_{k=1}^m \bar{u}_k Y_k \rangle} + O(||Z||^3 t^5) \\
&=& \sum_{i=1}^m \bar{u}_i Y_i - \frac{1}{2} \overline{\sum_{j=1}^m \bar{u}_j^2 \langle Y_j : Y_j \rangle + 2 \sum_{j<k} \bar{u}_j \bar{u}_k \langle Y_j:Y_k \rangle} + O(||Z||^3 t^5) \, ,
\end{eqnarray*}
where $Z=\sum_{i=1}^m \bar{u}_i Y_i$. Now, observe that $\displaystyle{\frac{d}{dt} (\phi(c(t)))}= \dot{c}(t) (\phi)$. Then, using properties {\it (iii)} and {\it (iv)} of $\phi$, we get
\begin{eqnarray*}
\frac{d}{dt} (\phi(c(t))) = \frac{1}{2} \overline{\sum_{j=1}^m a_{jj} \bar{u}_j^2 + 2 \sum_{j<k}a_{jk}\bar{u}_j \bar{u}_k} 
+ O(||Z||^3 t^5) \, .
\end{eqnarray*}
The expression $\sum_{j=1}^m a_{jj} \bar{u}_j^2 + 2 \sum_{j<k}a_{jk}\bar{u}_j \bar{u}_k$ does not change its sign, whatever the functions $u_1(t),\dots, u_m(t)$ might be, because as a quadratic polynomial in $\bar{u}_{s_1}$ its radicand is always negative. Therefore, $\displaystyle{\frac{d}{dt} (\phi(c(t)))}$ has constant sign for sufficiently small $t$, since $\overline{\sum_{j=1}^m a_{jj} \bar{u}_j^2 + 2 \sum_{j<k}a_{jk}\bar{u}_j \bar{u}_k} = O(t^3)$ and dominates $O(||Z||^3 t^5)$. Finally,
\[
\phi (c(t)) = \phi (q_0) + \int_0^t \frac{d}{ds} (\phi(c(s))) = \int_0^t \frac{d}{ds} 
(\phi(c(s)))
\]
will have constant sign for $t$ small enough. As a consequence, all the points 
in a neighbourhood of $q_0$ where $\phi$ has the opposite sign (property {\it (v)}) are 
unreachable in small time, which contradicts the hypothesis of controllability.

If $(a_{s_{m-1}s_m}^{(m-1)})^2 - a_{s_{m-1}s_{m-1}}^{(m-1)} a_{s_ms_m}^{(m-1)} = 0$, we can do 
the following. For $j=1$, we choose $b_{1s_m} \not = 0$ and
\begin{eqnarray}\label{descarte}
b_{1s_{m-1}} &=& -b_{1s_{m}} \frac{a^{(m-1)}_{s_{m-1}s_m}}{a^{(m-1)}_{s_{m-1}s_{m-1}}} = C_{s_{m-1}} b_{1s_m} \nonumber \\
b_{1s_{m-2}} &=& -\frac{a^{(m-2)}_{s_{m-2}s_{m-1}}b_{1s_{m-1}} + a^{(m-2)}_{s_{m-2}s_{m}} 
b_{1s_{m}}}{a^{(m-2)}_{s_{m-2}s_{m-2}}} = C_{s_{m-2}} b_{1s_m} \nonumber \\
& \vdots & \\
b_{1s_1}&=&-\frac{\sum_{k \not=s_1} b_{1k} a_{ks_1}}{a_{s_1s_1}} = C_{s_{1}} b_{1s_m} \, . \nonumber
\end{eqnarray}
We denote $C_{s_m}=1$. For $j > 1$, we select the $(b_{jk})_{1 \le k \le m}$ such that the matrix $B$ be non-singular. Consequently, we change our original basis $\{ Y_1,\dots, Y_m\}$ to a new one $\{ Y'_1,\dots, Y'_m\}$. In this basis, following (\ref{eluno}), one has
\begin{eqnarray*}
\langle Y_1':Y_1' \rangle (q_0)&=& lc(Y'_1(q_0),\dots,Y'_m(q_0)) \\
\langle Y_j':Y_j' \rangle (q_0)&=& lc(Y'_1(q_0),\dots,Y'_m(q_0)) + a_{jj}' \left \langle Y_1:Y_2 \right \rangle (q_0)\, , \; \; 2 \le j \le m \, .
\end{eqnarray*}
In addition, one can check that
\begin{eqnarray*}
\langle Y_1':Y'_j \rangle (q_0) &=& lc(Y'_1(q_0),\dots,Y'_m(q_0)) + \left(\sum_{k,l} a_{kl}b_{1k}b_{jl} \right) \langle Y_1:Y_2 \rangle (q_0) \\
&=& lc(Y'_1(q_0),\dots,Y'_m(q_0)) + b_{1s_m} \left( \sum_{l} b_{jl} \left( \sum_{k} a_{kl} C_k \right) \right) \langle Y_1:Y_2 \rangle (q_0) \, , \; 2 \le j \le m \, .
\end{eqnarray*}
Now, computations with {\it Mathematica} show that when the $C_{k}$ are given by (\ref{descarte}), then
\[
\sum_{k} a_{kl} C_k = 0 \, , \; \; 1 \le l \le m \, ,
\]
and this guarantees us that
\[
\langle Y_1':Y'_j \rangle (q_0) = lc(Y'_1(q_0),\dots,Y'_m(q_0)) \, , \; \; 2 \le j \le m \, .
\]
If the $a_{jj}' = 0$, $2 \le j \le m$, we are done. Assume then that $a_{33}' \not = 0$, reordering the input vector fields if necessary. Assume further that $\langle Y_2':Y_3' \rangle (q_0)$ is not a linear combination of $\{ Y'_1,\dots, Y'_m\}$ (otherwise, redefine a new $Y_2''$ as $Y_2'+Y_3'$). Then we have,
\begin{eqnarray*}
\left \langle Y_2':Y_2' \right \rangle
(q_0) &=& lc (Y_1'(q_0),\dots, Y_m'(q_0)) + a_{11}' \left \langle Y_2':Y_3'
\right \rangle (q_0) \\
& \vdots & \\
\left \langle Y_m':Y_m' \right \rangle
(q_0) &=& lc (Y_1'(q_0),\dots, Y_m'(q_0)) + a_{mm}' \left \langle Y_2':Y_3'
\right \rangle (q_0) \\
\left \langle Y_2':Y_4' \right \rangle (q_0) &=& lc
(Y_1'(q_0),\dots, Y_m'(q_0)) + a_{24}'\left \langle Y_2':Y_3' \right \rangle
(q_0) \\
& \vdots & \\
\left \langle Y_{m-1}':Y_m' \right \rangle (q_0) &=& lc
(Y_1'(q_0),\dots, Y_m'(q_0)) + a_{m-1m}' \left \langle Y_2':Y_3' \right \rangle
(q_0) \, ,
\end{eqnarray*}
where we have denoted with a little abuse of notation by $a_{jk}'$ the new coefficients corresponding to $\langle Y_2' : Y_3' \rangle$. Consequently, we can now reproduce the preceding discussion, but with the $m-1$ vector fields $\{ Y_2', \dots , Y_m' \}$, since $Y_1'$ does not affect the situation. We look for one change of basis $B'$ in the vector fields $\{ Y_2', \dots, Y_m' \}$ such that the new ones $\{ Y_2'', \dots, Y_m'' \}$ together with $Y_1'$ verify the sufficient conditions for STLCC. Accordingly, we must consider the vanishing of the polynomials
\[
\sum_{k=2}^{m} {b_{jk}^2}' a_{kk}' + 2 \hspace{-8pt} \sum_{2 \le k < l \le m} b_{jk}' b_{jl}' a_{kl}' = 0 \, , \; \; 2 \le j \le m \, .
\]
The cases in which the last radicand $({a_{s_{m-1}s_m}^{(m-1)}}')^2 - {a_{s_{m-1}s_{m-1}}^{(m-1)}}' {a_{s_ms_m}^{(m-1)}}'$ does not vanish are treated as before. When it vanishes, we obtain a new basis $\{ Y_1''=Y_1', Y_2'', \dots, Y_m''\}$ such that
\begin{eqnarray*}
\langle Y_1'':Y_1'' \rangle (q_0) \hspace{-10pt} &,& \hspace{-7pt}\langle Y_2'':Y_2'' \rangle (q_0) = lc(Y''_1(q_0),\dots,Y''_m(q_0)) \\
\langle Y_j'':Y_j'' \rangle (q_0)&=& lc(Y''_1(q_0),\dots,Y''_m(q_0)) + c_{jj}' \left \langle Y_2':Y_3' \right \rangle (q_0)\, , \; \; 3 \le j \le m \\
\langle Y_1'':Y''_j \rangle \hspace{-10pt} &,& \hspace{-7pt} \langle Y_2'':Y''_{j+1} \rangle = lc(Y''_1(q_0),\dots,Y''_m(q_0)) \, , \; \; 2 \le j \le m \, ,
\end{eqnarray*}
where possibly there exits some $3 \le j \le m$ such that $c_{jj}' \not = 0$. By an induction procedure, we finally come to consider discarding the case of a certain basis 
$\{Z_1=Y_1',Z_2=Y_2'',\dots,Z_m\}$ of $D$ satisfying $\langle Z_i:Z_j \rangle (q_0) 
\in \hbox{span} \{Z_1 (q_0),\dots , Z_m(q_0)\}$, $1 \le i < j \le m$, and 
the sufficient conditions for STLCC at $q_0$ for $Z_1, \dots, Z_{m-1}$, but 
such that $\langle Z_m : Z_m \rangle (q_0) \not \in \hbox{span} \{Z_1 (q_0), 
\dots , Z_m(q_0) \}$. Similarly as we have done above, the application of Lemma 
\ref{magic} with the vector fields $\{ Z_1 , \dots, Z_m, \langle Z_m:Z_m \rangle\}$ implies that the system is not controllable at $q_0$, yielding a contradiction.

The preceding discussion has been made on the assumption that a series of non-zero elements $(a^{(i)}_{s_is_i})_{1 \le i \le m-1}$ exists. It remains to consider the possibility when there exists an $i \ge 2$ such that $a^{(i)}_{kk} = 0$, for all $k \in \{ 1,\dots,m \} / \{ s_1, \dots, s_{i-1} \}$. This means that the polynomial
\[
\sum_{k \not = s_1,\dots,s_{i-1}} b_{jk}^2 a_{kk}^{(i)} + 2 \hspace*{-23pt} \sum_{k<l,k,l \not = s_1,\dots,s_{i-1}} \hspace*{-17pt} b_{jk} b_{jl} a_{kl}^{(i)} 
\]
is no longer quadratic, but it has the form
\begin{equation}\label{nextII}
2 \hspace*{-23pt} \sum_{k<l,k,l \not = s_1,\dots,s_{i-1}} \hspace*{-17pt} b_{jk} b_{jl} a_{kl}^{(i)} \, .
\end{equation}
If any of the $a_{kl}^{(i)}$ is different from zero, then it is clear that we can choose the 
$b_{jk}$, $k \not \in \{ s_1,\dots,s_{i-1} \}$, such that (\ref{nextII}) be positive. Then, reasoning as before, we find a regular matrix $B$ yielding the desired change of basis. If this is not the case, i.e. $a_{kl}^{(i)} = 0$, for all $k<l,k,l \not \in \{ s_1,\dots,s_{i-1} \}$, we can do the following. Choose $\{ (b_{jk})_{1\le j \le m}\}$, with $k \not \in \{ s_1,\dots,s_{i-1} \}$ $m-i+1$ linearly independent vectors in $\R^m$ such that the minor $\{ b_{jk} \}_{1\le j \le m-i+1 }^{k \not = s_1, \dots, s_{i-1}}$ is regular. Now, let in (\ref{descarte}) that $j$ varies between $1$ and $m-i+1$, that is, take
\begin{eqnarray}\label{descarteII}
b_{js_{i-1}} &=& -\frac{\sum_{k \not = s_1 ,\dots, s_{i-1}}^{m} b_{jk}
a^{(i-1)}_{s_{i-1}k}}{a^{(i-1)}_{s_{i-1}s_{i-1}}} \nonumber \\
b_{js_{i-2}} &=& -\frac{\sum_{k \not = s_1 ,\dots, s_{i-2}}^{m} b_{jk}
a^{(i-2)}_{s_{i-2}k}}{a^{(i-2)}_{s_{i-2}s_{i-2}}} \nonumber \\ & \vdots
& \\ b_{js_1}&=&-\frac{\sum_{k \not =s_1} b_{jk} a_{ks_1}} {a_{s_1s_1}} \, , \nonumber
\end{eqnarray}
for $1 \le j \le m-i+1$. Finally, for $j > m-i+1$, we select the $b_{jk}$ such that the matrix $B$ is non-singular. In this manner, in an unique step, we would change to a new basis $\{ Y'_1,\dots, Y'_m\}$
verifying
\begin{eqnarray*}
\langle Y_1':Y_1' \rangle (q_0) \hspace{-10pt} &,& \hspace{-7pt} \dots \, , \langle Y_{m-i+1}':Y_{m-i+1}' \rangle (q_0) =
lc(Y'_1(q_0),\dots,Y'_m(q_0)) \\
\langle Y_j':Y_j' \rangle (q_0)&=& lc(Y'_1(q_0),\dots,Y'_m(q_0)) + a_{jj}' \left \langle Y_1:Y_2 \right \rangle (q_0)\, , \; \; m-i+1 \le j \le m \\
\langle Y_k':Y_l' \rangle (q_0)&=& lc(Y'_1(q_0),\dots,Y'_m(q_0)) \, , \; \; k < l, 1 \le k \le m-i+1 \, ,
\end{eqnarray*}
with possibly some of the $(a_{jj}')_{m-i+1 \le j \le m}$ being different from zero. Now,
the above discussion can be redone in this context to assert the validity of the theorem.

In short, what we have done is the following: firstly, we have considered the case when for all $i$ there exists a $s_i$ such that $a_{s_is_i}^{(i)} \not = 0$. We have seen that this case can be subdivided into three: one ensuring the desired change of basis, another one which is not possible under the hypothesis of small-time local configuration controllability and a third one is a kind of ``reduced" situation where we can get rid of the problems caused by one vector field. Then, under the same assumption on the new coefficients, $a_{jk}'$ (i.e. for all $i$, there exists a $s_i$ such that ${a^{(i)}}'_{s_is_i} \not = 0$), we can reproduce the discussion. If we repeteadly fall into the third case, we eventually find a contradiction with the controllability assumption. Secondly, we have treated the case when this type of ``circular" process is broken: that is, there exists an $i$ such that $a_{kk}^{(i)}= 0$, for all $k \not = s_1,\dots,s_{i-1}$. What we have shown then is that this leads to the obtention of either new input vector fields satisfying the sufficient conditions for STLCC or a reduced situation where we can get rid at the same time of the problems associated to $m-i+1$ vector fields.
\QED

\begin{remark}
{\rm Notice that the proof of this result can be reproduced for the corresponding notions of accessibility and controllability at zero velocity. Indeed, a mechanical control system of the form (\ref{affine}) with $m=n-1$, which is STLC at $q_0$ and zero velocity is in particular STLCC at $q_0$. Then, Theorem \ref{ole} implies that there exists a basis of input vector fields ${\cal Y}$ satisfying the sufficient conditions of Theorem \ref{suffi}, so the same result is also valid for local controllability at zero velocity.
}
\end{remark}

\begin{remark}
{\rm
We would like to point out that the differential flatness properties of this type of underactuated mechanical control systems have also been studied in intrinsic geometric terms. Indeed, in \cite{RaMu} it was proved a characterization of configuration flatness for Lagrangian control systems with $n$ degrees of freedom and $n-1$ controls and without nonholonomic constraints. The importance of flatness to control applications lies in the fact that it provides a systematic and simple way to generate trajectories between two given states.
}
\end{remark}

\section{Open problem}

So far, we have not succeeded in generalizing the proof of Theorem \ref{ole} to the full general case, that is, a system with $n$ degrees of freedom and $m <n$ input forces or vector fields. In fact, the following simplest case one can consider presents already serious difficulties. Let (\ref{affine}) be a mechanical control system on a 4-dimensional configuration manifold with $m=2$ inputs. Assume it is STLCC at $q_0 \in Q$. Let $Y_1$, $Y_2 \in D$ be linearly independent vector fields such that $ \hbox{span} \{ Y_1(q_0), Y_2(q_0), \langle Y_1 : Y_1 \rangle (q_0) ,\langle Y_1 : Y_2  \rangle (q_0)\} = T_{q_0} Q$. Can one discard the possibility
\[
\langle Y_2 : Y_2  \rangle (q_0) = a_1 Y_1 (q_0) + a_2 Y_2 (q_0) + a_3 \langle Y_1 : Y_2  \rangle (q_0)+ a_4 \langle Y_1 : Y_1  \rangle (q_0) \, ,
\]
with $a_4 <0$ and $a_3^2 + 4a_4 <0$?

The conditions on the coefficients $a_3$, $a_4$ prevent us from finding an adequate change of basis of $D_{q_0}$ such that the new vector fields satisfy the sufficient conditions for STLCC, or simply that make the situation clearer. Moreover, Lemma \ref{magic} is of no help in this case since when applying it to the vector fields $\{ Y_1, Y_2, \langle Y_1:Y_2 \rangle , \langle Y_2:Y_2 \rangle\}$ we get
\begin{eqnarray*}
\dot{c}(t) (\phi) &=& \left( \bar{u}_1 Y_1 + \bar{u}_2
Y_2 - \frac{1}{2}\overline{(\bar{u}_1^2 \left \langle Y_1:Y_1\right \rangle + 2 \bar{u}_1\bar{u}_2 \left \langle Y_1:Y_2\right \rangle + \bar{u}_2^2 \left \langle Y_2:Y_2\right \rangle)} + O(\|Z\|^3 t^5) \right) (\phi) \\
&=& \frac{1}{2} \overline{(\bar{u}_1^2 + a_4 \bar{u}_2^2)} + O(\|Z\|^3 t^5) \, ,
\end{eqnarray*}
where $Z=u_1 Y_1 + u_2 Y_2$. Clearly, $\bar{u}_1^2 + a_4 \bar{u}_2^2$ can change its sign. Another possibility is to try to sharpen Lemma \ref{magic} with a view in this case. Without much more difficulty, one can modify the proof of the Lemma to ensure, under the assumptions given, the existence of a function $\phi$ verifying properties {\it (i), (ii), (v)} and

{\it (iii)}' $Y_1 (\phi) = Y_2 (\phi) = 0$ on a neighbourhood $V$ of $q_0$

{\it (iv)}' $\langle Y_1 : Y_2 \rangle (\phi) (q_0) = \pm 1$, $\langle Y_1 : Y_1 \rangle (\phi) (q_0) = -1$.

But again this shares the same lot, since we obtain
\begin{eqnarray*}
\dot{c}(t) (\phi) &=& \frac{1}{2} \overline{\left( \bar{u}_1^2 \mp 2 \bar{u}_1 \bar{u}_2 + \bar{u}_2^2 (\mp a_3 + a_4) \right)} \\
&=& \frac{1}{2} \overline{ (\bar{u}_1 - \bar{u}_2 (\pm 1 + \sqrt{1-a_4 \pm a_3})) (\bar{u}_1 - \bar{u}_2 (\pm 1 - \sqrt{1-a_4 \pm a_3}))}\, .
\end{eqnarray*}
This expression can change its sign due to $a_3^2 + 4 a_4 <0 \Rightarrow |a_3| < 1-a_4$.

The following simple example fits perfectly in the above exposed case, except for that we do not know if it is small time configuration controllable. If it were so, it would constitute a counter example of Lewis' conjecture in its general formulation.

\begin{example}
{\rm
Consider a mechanical control system on $Q=\R^4$, with coordinates $(x,y,z,w)$. The Riemannian metric is given by
\[
g=dx \otimes dx + dy \otimes dy + dz \otimes dz + dw \otimes dw \, ,
\]
and the input vector fields
\begin{eqnarray*}
Y_1 &=& (1+z) \frac{\partial}{\partial x} + \frac{\partial}{\partial y} + \frac{\partial}{\partial z}+ (1+y) \frac{\partial}{\partial w} \, , \\
Y_2 &=& \frac{\partial}{\partial y} - 2 \frac{\partial}{\partial z} - (1+y) \frac{\partial}{\partial w} \, .
\end{eqnarray*}
In coordinates, the control equations are given by
\begin{equation}\label{mate}
\left\{
\begin{array}{rcl}
\ddot{x} &=& (1+z) u_1 \\
\ddot{y} &=& u_1 + u_2 \\
\ddot{z} &=& u_1 - 2 u_2\\
\ddot{w} &=& (1+y) (u_1 - u_2) \, .
\end{array}
\right.
\end{equation}
Since
\begin{eqnarray*}
\left \langle Y_1 : Y_1 \right \rangle &=& 2 \left( \frac{\partial}{\partial x} +  \frac{\partial}{\partial w} \right) \\
\left \langle Y_1 : Y_2 \right \rangle &=& -2 \frac{\partial}{\partial x} \\
\left \langle Y_2 : Y_2 \right \rangle &=& -2 \frac{\partial}{\partial w} \, ,
\end{eqnarray*}
we deduce that $\hbox{span} \{Y_1(q), Y_2(q), \langle Y_1:Y_2 \rangle (q), \langle Y_1:Y_1 \rangle (q) \} = T_qQ$ for all $q \in Q$ and the system (\ref{mate}) is locally configuration accessible.

In addition, $\langle Y_2 : Y_2 \rangle = - \langle Y_1 : Y_2 \rangle - \langle Y_1 : Y_1 \rangle$, so $a_4= -1 <0$ and $a_3^2 +4 a_4=-3<0$ and we can not conclude that (\ref{mate}) is locally configuration controllable. In this respect, observe that whatever projection of the control system (\ref{mate}) to 3 dimensions we consider, that is, in the variables $(x,y,z)$ or $(y,z,w)$, we obtain that the projected system is STLCC. Indeed, we can find a change of input vector fields such that the new ones satisfy the sufficient conditions (Theorem \ref{ole}, $n-1=3$, $m=2$).
}
\end{example}

\section{Conclusions}

In this paper, we have been able to demonstrate that the sufficient conditions encountered in \cite{LeMu2} for a mechanical control system to be STLCC are also necessary when the configuration manifold is $n$-dimensional and there are $n-1$ inputs, in the sense that there exists some basis of input vector fields that verifies them.

In fact, the algorithmic nature of the proof of Theorem \ref{ole} allows us to find such a basis. Perhaps it would be useful to determine an ordered procedure that, for a given mechanical system, leads us to the obtention of such a basis in a systematic way. This could be very interesting in a number of applications, including motion planning, trajectory tracking, reliability of autonomous vehicles and the design of mechanisms with fewer actuators than typical, yielding less costly devices.

On the other hand, as we point out in Section 6, the validity of the conjecture in the full general case remains still open.

\section*{Acknowledgements}
This work was partially supported by FPU and FPI grants from the Spanish Ministerio de Educaci\'on y Cultura and grant DGICYT PB97-1257. We wish to thank F. Cantrijn for several useful suggestions and the Department of Mathematical Physics and Astronomy of the University of Ghent for its kind hospitality.

\bigskip

{\parindent 0cm

{\sc Jorge Cort{\'e}s \dag, Sonia Mart{\'\i}nez \ddag} 

{\it Laboratory of Dynamical Systems, Mechanics and Control \\
Instituto de Matem\'aticas y F{\'\i}sica Fundamental, CSIC \\
Serrano 123, 28006 Madrid, SPAIN\\
\dag e-mail: j.cortes@imaff.cfmac.csic.es \quad 
\ddag e-mail: s.martinez@imaff.cfmac.csic.es \\}

\end{document}